\documentclass{amsart}

\usepackage{amssymb}
\newtheorem{thm}{Theorem}
\newtheorem{lem}{Lemma}

\begin{document}

\bibliographystyle{plain}

\title[Romi Shamoyan ]{On some extremal problems in certain  harmonic function spaces of several variables related to mixed norm spaces}

\author[]{Romi Shamoyan}

\address{Department of Mathematics, Bryansk State Technical University, Bryansk 241050, Russia}
\email{\rm rshamoyan@gmail.com}

\date{}

\begin{abstract}
Some  estimates on distances in spaces of harmonic functions in the unit ball and the upper half space are provided.
New estimates concerning mixed norm spaces and general weighted Bergman spaces are obtained  and discussed.

\end{abstract}

\maketitle

\footnotetext[2]{\, Mathematics Subject Classification 2010 Primary 42B15, Secondary 42B30.  Key words
and Phrases: Distance estimates,harmonic function,unit ball,Bergman spaces}

\section{Introduction and preliminaries}

The main goal of this note to  present several new results on distances in harmonic spaces though some isolated results from our previous papers (see\cite{AS1},\cite{AS2})
will be also given to make the picture more complete.
This direction of investigation related to 
extremal problems in harmonic function spaces started in \cite{AS1},then continued in \cite{AS2},\cite{AS3}
In this note in  particular a new  sharp theorem on distances in harmonic function spaces in the unit ball will be added.An analogue of this assertion in
case of upper half space was formulated and proved before in \cite{AS2}
and for completeness of exposition we will add it also here (but without proof)
Note also this  result in the unit ball in case of $p>1$  was also provided before in \cite {AS2} ,in this note we extend it  to all values of  positive $p$ .
The main tool for  us will be  so-called Whitney type decomposition of the unit ball of Euclidean space and this technique was used before
in case of upperhalfspace  by us in \cite{AS2} 
We also add important remarks concerning this extremal problem  (distance function) in mixed norm spaces and weighted spaces of harmonic functions at the end of this paper.
Harmonic function spaces were studied by many authors during last several decades we mention , for example, \cite{JP1} and \cite{DS} and various references there.

Our line of investigation can be also considered as a continuation of papers on distances in analytic function spaces \cite{SM1} and \cite{SM2}. All main results are 
contained in the second section of the paper. The first section is devoted to preliminaries and main definitions which are needed for formulations of main results.
Almost all objects we define in this section and definitions can be found in \cite{DS},\cite {AS1},\cite{AS2},\cite{AS3} and in \cite{St}.The Whitney decomposition and properities and estimates we used for our proofs 
can be found in various places, but we refer the reader to \cite {AS1} ,\cite {AV} for case of upperhalfspace and for case of unit ball of Euclidean space (see also \cite{AS2})
This paper also contains  important additions namely we add important remarks about general Bergman spaces with general $w$ weights,these new weighted harmonic spaces appeared for the first time in \cite{MZ}) .These results are also new and 
we hope to return to these issues related to general harmonic  $h^p_w$ and $H^{p}_w$ Bergman- type spaces in unit ball and upper halfspace also later.  

Let $\mathbb B$ be the open unit ball in $\mathbb R^n$, $\mathbb S = \partial \mathbb B$ is the unit sphere in $\mathbb R^n$, for $x \in \mathbb R^n$ we have $x = rx'$, where $r = |x| = \sqrt{\sum_{j=1}^n x_j^2}$ and $x' \in \mathbb S$. Normalized Lebesgue measure on $\mathbb B$ is denoted by $dx = dx_1 \ldots dx_n = r^{n-1}dr dx'$ so that $\int_{\mathbb B} dx = 1$. We denote the space of all harmonic functions in an open set $\Omega$ by $h(\Omega)$. In this paper letter $C$ designates a positive constant which can change its value even in the same chain of inequalities.

For $0<p<\infty$, $0 \leq r < 1$ and $f \in h(\mathbb B)$ we set
$$M_p(f, r) = \left( \int_{\mathbb S} |f(rx')|^p dx' \right)^{1/p},$$
with the usual modification to cover the case $p = \infty$. 

For $0<p<\infty$ and $\alpha > -1$ we consider weighted harmonic Bergman spaces $A^p_\alpha = A^p_\alpha(\mathbb B)$ defined by
$$A^p_\alpha = \left\{ f \in h(\mathbb B) : \| f \|_{A^p_\alpha}^p = \int_{\mathbb B} |f(x)|^p (1-|x|^2)^\alpha dx < \infty
\right\}$$
For $p = \infty$ this definition is modified in a standard manner: 
$$A^\infty_\alpha = A^\infty_\alpha(\mathbb B) = \left\{ f \in h(\mathbb B) : \| f \|_{A^\infty_\alpha} = \sup_{x \in \mathbb B} |f(x)| (1-|x|^2)^\alpha < \infty \right\}, \quad \alpha > 0$$

These spaces are complete metric spaces for $0<p\leq \infty$, they are Banach spaces for $p \geq 1$ (\cite{DS})
These spaces serve as particular cases of more general scales of mixed norm spaces $F^{p,q}_\alpha$ and $B^{p,q}_\alpha$ (see \cite {AS2})in unit ball of $R^{n}$

Next we need certain facts on spherical harmonics and the Poisson kernel, see \cite{DS} for a detailed exposition. Let $Y^{(k)}_j$ be the spherical harmonics of order $k$, $1 \leq j \leq d_k$, on $\mathbb S$. Next,
$$Z_{x'}^{(k)}(y') = \sum_{j=1}^{d_k} Y_j^{(k)}(x') \overline{Y_j^{(k)}(y')}$$
are zonal harmonics of order $k$. Note that the spherical harmonics $Y^{(k)}_j$, ($k \geq 0$, $1 \leq j \leq d_k$) form an orthonormal basis of $L^2(\mathbb S, dx')$. Every $f \in h(\mathbb B)$ has an expansion
$$f(x) = f(rx') = \sum_{k=0}^\infty r^k b_k\cdot Y^k(x'),$$
where $b_k = (b_k^1, \ldots, b_k^{d_k})$, $Y^k = (Y_1^{(k)}, \ldots, Y_{d_k}^{(k)})$ and $b_k\cdot Y^k$ is interpreted in the scalar product sense: $b_k \cdot Y^k = \sum_{j=1}^{d_k} b_k^j Y_j^{(k)}$.

We denote the Poisson kernel for the unit ball as usual by $P(x, y')$, it is given by
\begin{align*}
P(x, y') = P_{y'}(x) & = \sum_{k=0}^\infty r^k \sum_{j=1}^{d_k} Y^{(k)}_j(y') Y^{(k)}_j(x') \\
& = \frac{1}{n\omega_n} \frac{1-|x|^2}{|x-y'|^n}, \qquad x = rx' \in \mathbb B, \quad y' \in \mathbb S,
\end{align*}
where $\omega_n$ is the volume of the unit ball in $\mathbb R^n$. We are going to use also a Bergman kernel for $A^p_\beta$ spaces,(see \cite{DS}) this is the following function
\begin{equation}\label{bker}
Q_\beta(x, y) = 2 \sum_{k=0}^\infty \frac{\Gamma(\beta + 1 + k + n/2)}{\Gamma(\beta + 1) \Gamma(k + n/2)}
r^k \rho^k Z_{x'}^{(k)}(y'), \qquad x = rx', \; y = \rho y' \in \mathbb B.
\end{equation}
This function is playing a very important role in many issues related with harmonic function spaces in unit ball (see \cite{DS},\cite{AS1},\cite{AS2},\cite{AS3} and references there

For details on this kernel we refer to \cite{DS}, where the following theorem can be found.
\begin{thm}[\cite{DS}]\label{intrep}
Let $p \geq 1$ and $\beta \geq 0$. Then for every $f \in A^p_\beta$ and $x \in \mathbb B$ we have
$$f(x) = \int_0^1 \int_{\mathbb S^{n-1}} Q_\beta(x, y) f(\rho y')(1-\rho^2)^\beta \rho^{n-1}d\rho dy',
\qquad y = \rho y'.$$
\end{thm}

This theorem is a cornerstone for our approach to distance problems in the case of the unit ball. 
Everywhere below $y=\rho y"$ and $x=rx"$ are in unit ball, sometimes this will be omitted by us and this will be clear from text.
The following lemma  gives estimates for this kernel, see \cite{DS}, \cite{JP1}.
\begin{lem}\label{DSlemma}
1. Let $\beta > 0$. Then, for $x = rx', y = \rho y' \in \mathbb B$ we have
$$|Q_\beta(x, y)| \leq \frac{C}{|\rho x - y'|^{n+\beta}}.$$
2. Let $\beta > -1$. Then
$$\int_{\mathbb S^{n-1}} |Q_\beta(rx', y)| dx' \leq \frac{C}{(1-r\rho)^{1+\beta}},
\qquad |y| = \rho, \quad 0 \leq r < 1.$$
3. Let $\beta > n-1$, , $0 \leq r < 1$ and $y'\in \mathbb S^{n-1}$. Then
$$\int_{\mathbb S^{n-1}} \frac{dx'}{|rx' - y'|^\beta} \leq \frac{C}{(1-r)^{\beta-n+1}}.$$
\end{lem}

The following simple lemma is classical see ,for example, \cite {DS} and references there.

\begin{lem}[\cite{DS}]\label{rro}
Let $\alpha > -1$ and $\lambda > \alpha + 1$. Then
$$\int_0^1 \frac{(1-r)^\alpha}{(1-r\rho)^\lambda} dr \leq C (1-\rho)^{\alpha + 1 - \lambda}, \qquad 0 \leq \rho < 1.$$
\end{lem}

The following lemma is purely technical ,but it is vital for the proof of our main result.

\begin{lem}\label{qbeta}
For $\delta > -1$, $\gamma > n + \delta$ and $\beta > 0$ we have
$$\int_{\mathbb B} |Q_\beta(x, y)|^{\frac{\gamma}{n+\beta}} (1-|y|)^\delta dy \leq C (1-|x|)^{\delta - \gamma + n},
\qquad x \in \mathbb B.$$
\end{lem}

{\it Proof.} Using Lemma \ref{DSlemma} and Lemma \ref{rro} we obtain:
\begin{align*}
\int_{\mathbb B} |Q_\beta(x, y)|^{\frac{\gamma}{n+\beta}} (1-|y|)^\delta dy & \leq C \int_{\mathbb B}
\frac{(1-|y|)^\delta}{|\rho rx' - y'|^\gamma} dy \\
& \leq C \int_0^1 (1-\rho)^\delta \int_{\mathbb S} \frac{dy'}{|\rho r x' - y' |^\gamma} dy' d\rho \\
& \leq C \int_0^1 (1-\rho)^\delta (1-r\rho)^{n - \gamma - 1}d\rho \\
& \leq C (1-r)^{n + \delta - \gamma}. \qquad \Box
\end{align*}

Now we turn to the basic definitions  for upperhalfspace,all of them are classical see \cite{DS} and references there or see \cite{AS1},\cite{AS2}

We first set $\mathbb R^{n+1}_+ = \{ (x, t) : x \in \mathbb R^n, t > 0 \} \subset \mathbb R^{n+1}$. We usually denote points in $\mathbb R^{n+1}_+$ by $z = (x, t)$ or $w = (y, s)$ where $x, y \in \mathbb R^n$ and $s, t > 0$.

For $0 < p < \infty$ and $\alpha > -1$ we consider spaces
$$\tilde A^p_\alpha(\mathbb R^{n+1}_+) = \tilde A^p_\alpha = \left\{ f \in h(\mathbb R_+^{n+1}) : \int_{\mathbb R_+^{n+1}} |f(x, t)|^p t^\alpha dx dt < \infty \right\}.$$
Also, for $p = \infty$ and $\alpha > 0$, we set
$$\tilde A^\infty_\alpha(\mathbb R^{n+1}_+) = \tilde A^\infty_\alpha = \left\{ f \in h(\mathbb R_+^{n+1}) :
\sup_{(x, t) \in \mathbb R^{n+1}_+} |f(x, t)|t^\alpha < \infty \right\}.$$
These spaces have natural (quasi)-norms, for $1 \leq p \leq \infty$ they are Banach spaces and for $0<p\leq 1$ they are
complete metric spaces.(see \cite{DS})
They serve as paticular cases of more general $F^{p,q}_{\alpha}$ and $B^{p,q}_\alpha$ mixed norm spaces. (see ,for example, \cite{AS2} and references there)

Now we turn to formualte some known assertions in this case of upperhalfspaces (see \cite{DS}) and references there
We denote the Poisson kernel for $\mathbb R^{n+1}_+$ by $P(x, t)$, i.e.
$$P(x, t) = c_n \frac{t}{(|x|^2 + t^2)^{\frac{n+1}{2}}}, \qquad x \in \mathbb R^n, t > 0.$$
For an integer $m \geq 0$ we introduce a Bergman kernel $Q_m(z, w)$, where $z = (x, t) \in \mathbb R^{n+1}_+$ and $w = (y, s) \in \mathbb R^{n+1}_+$, by
$$Q_m(z, w) = \frac{(-2)^{m+1}}{m!} \frac{\partial^{m+1}}{\partial t^{m+1}} P(x-y, t+s).$$

The terminology is justified by the following result from \cite{DS} which is a complete analogue 
of integral representation of Bergman spaces in the unit ball we formulated above.Note it is well-known that these theorems 
on integral representations for ball and upperhalfspace have various applications in harmonic function theory. 

\begin{thm}\label{brthm}
Let $0<p<\infty$ and $\alpha > -1$. If $0<p\leq 1$ and $m \geq \frac{\alpha +n +1}{p} - (n+1)$ or $1\leq p < \infty$
and $m> \frac{\alpha +1}{p} - 1$, then
\begin{equation}\label{brep}
f(z) = \int_{\mathbb R^{n+1}_+} f(w) Q_m(z, w) s^m dy ds, \qquad f \in \tilde A^p_\alpha,\quad z \in \mathbb R^{n+1}_+.
\end{equation}
\end{thm}

The following elementary estimate of this kernel is crucial and it is contained, for example, in \cite{DS}:
\begin{equation}\label{estq}
|Q_m(z, w)| \leq C \left[ |x-y|^2 + (s+t)^2 \right]^{-\frac{n+m+1}{2}}, \quad z = (x, t), w = (y, s) \in \mathbb R^{n+1}_+.
\end{equation}

It is well- known that the theory of Bergman spaces in unit ball and upperhalfspace are parallel to each other and it is natural to consider extremal problems in both spaces together.
Various issues in these spaces can be solved using arguments related so- called properties of a certain family of special cubes which usually called Whitney cubes (see \cite{St}) and references there.
Here are two basic properties of Whitney type decomposition which we will need for upperhalfspace $R^{n+1}_{+}$ (see,for example, \cite{AS2}).
There is a collection of closed cubes $\Delta_{k}$ in $R^{n+1}_{+}$ with sides parallel to coordinate axes such that the following properties hold
The union of these collection $\Delta_{k}$ of cubes gives all $R^{n+1}_{+}$.The interior of $\Delta_{k}$ cubes are pairwise disjoint.Such type of family also
exists in the unit ball (see, for example, \cite {AV} and references there).Related facts and estimates we need for these cubes and their centers can be found in \cite{AV} in unit ball
and in \cite{AS2} in case of upperhalfspace  $R^{n+1}_{+}$

\section{Sharp estimates for distances in harmonic Bergman function spaces of several variables in the unit ball and in
$\mathbb R^{n+1}_+$ and related theorems in mixed norm spaces and weighted Bergman spaces.}

In this section we investigate distance problems both in the case of the unit ball and in the case of the upper
half space for harmonic functions of several variables. The method we use here originated in \cite{Zh}, see also \cite{AS1},\cite{AS2},\cite{AS3}, \cite{SM1}, \cite{SM2} for various modifications of this interesting approach,
and applications of this method to various distance problems in analytic function spaces and harmonic function spaces in one and several variables.
New assertions and important remarks concerning mixed norm harmonic function spaces and general harmonic weighted Bergman spaces will be added at the second part of this paper.
The followiing lemma is elementar.

\begin{lem}\label{emb}
Let $0 < p < \infty$ and $\alpha > -1$. Then there is a $C = C_{p, \alpha, n}$ such that for every $f \in A^p_\alpha
(\mathbb B)$ we have
$$|f(x)| \leq C (1-|x|)^{-\frac{\alpha + n}{p}} \| f \|_{A^p_\alpha}, \qquad x \in \mathbb B.$$
\end{lem}

{\it Proof.} We use subharmonic behavior of $|f|^p$, see for this for example \cite{St},\cite{DS}, to obtain
\begin{align*}
|f(x)|^p & \leq \frac{C}{(1-|x|)^n} \int_{B(x, \frac{1-|x|}{2})} |f(y)|^p dy\\
& \leq C \frac{(1-|x|)^{-\alpha}}{(1-|x|)^n} \int_{B(x, \frac{1-|x|}{2})} |f(y)|^p (1-|y|)^\alpha dy\\
& \leq C(1-|x|)^{-\alpha - n} \| f \|_{A^p_\alpha}^p. \qquad \Box
\end{align*}

This lemma shows that $A^p_\alpha$ is continuously embedded in $A^\infty_{\frac{\alpha + n}{p}}$ and motivates the
distance problem in unit ball which investigated in our Theorems  below.

\begin{lem}\label{lemb1}
Let $0<p<\infty$ and $\alpha > -1$. Then there is $C = C_{p, \alpha, n}$ such that for every $f \in \tilde A^p_\alpha$ and every $(x, t) \in \mathbb R^{n+1}_+$ we have
\begin{equation}\label{eqemb1}
|f(x, t)| \leq C y^{-\frac{\alpha + n + 1}{p}} \| f \|_{\tilde A^p_\alpha}.
\end{equation}
\end{lem}

The above lemma states that $\tilde A^p_\alpha$ is continuously embedded in $\tilde A^\infty_{\frac{\alpha + n +1}{p}}$, its proof is analogous to that of Lemma \ref{emb}.
This lemma motivates the distance problem in harmonic function spaces in upperhalfspaces which investigated  in our Theorems  below.

For $\epsilon > 0$, $t > 0$ and $f \in h(\mathbb B)$ we set
$$U_{\epsilon, t}(f) = U_{\epsilon, t} =  \{ x \in \mathbb B : |f(x)|(1-|x|)^t \geq \epsilon \}.$$
The following assertion alone can be found among other things in \cite{AS3} hovewer we decided  to formulate it since the remaining part of it $p\leq 1$ case was open and we
in this note close it in our next theorem,we also prove the first theorem to  use it for further results and observations at the end of this note.

\begin{thm}
Let $p > 1$, $\alpha > -1$, $\lambda = \frac{\alpha + n }{p}$, $m \in \mathbb N_0$ 
Set,
for $f \in  A^\infty_{\frac{\alpha + n }{p}}(\mathbb B)$:
$$s_1(f) = {\rm dist}_{ A^\infty_{\frac{\alpha + n }{p}}} (f,  A^p_\alpha),$$
$$s_2(f) = \inf \left\{ \epsilon > 0 : \int_{\mathbb B} \left( \int_{U_{\epsilon, \lambda}}
Q_\beta(x, y) (1-|y|)^{\beta - \lambda}dy \right)^p (1-|x|)^{\alpha} dx  < \infty \right\}.$$
Then there is a $m_{0}$ depending from $\alpha,p,n$ so that for all  $\beta > m_{0}$
 $$t_1(f) \asymp t_2(f)$$.
\end{thm}

In the following new sharp  theorem on distances we cover the remaining case of $ p\leq 1 $. Note this result is new.
Nevertheless the proof of this theorem is based on estimates connected with Whitney- type decomposition of the unit ball which we mentioned above(see \cite{AV}),and it is very close and parallel
to the case of upperhalfspace case (see the last theorem of this paper below) and arguments in proofs are  similar .
The  theorem in upperhalf space is formulated at the end of this paper and the proof of it is given in our previous work.
(see \cite{AS2}).Hence here we omit details of proof , referiyng the reader to \cite{AS2}.

\begin{thm}
Let $p \leq 1$, $\alpha > -1$, $\lambda = \frac{\alpha + n}{p}$, $m \in \mathbb N_0$ .
Set,
for $f \in A^\infty_{\frac{\alpha + n }{p}}(\mathbb B)$:
$$s_1(f) = {\rm dist}_{ A^\infty_{\frac{\alpha + n }{p}}} (f, A^p_\alpha),$$
$$s_2(f) = \inf \left\{ \epsilon > 0 : \int_{\mathbb B} \left( \int_{U_{\epsilon, \lambda}}
Q_\beta(x,y)(1-|y|)^{\beta-\lambda}dy \right)^p (1-|x|)^{\alpha} dx < \infty \right\}.$$
Then there is an $m_{0}$ depending on $\alpha,p,n$ so that  for all $\beta > m_{0}$
$$s_1(f) \asymp s_2(f)$$.
\end{thm}

We now provide now a complete proof of first theorem we formulated above.

{\it Proof.} We begin with inequality $t_1(f) \geq t_2(f)$. Assume $t_1(f) < t_2(f)$. Then there are
$0< \epsilon_1 < \epsilon$ and $f_1 \in A^p_\alpha$ such that $\|f-f_1\|_{A^\infty_t} \leq \epsilon_1$ and
$$\int_{\mathbb B} \left( \int_{U_{\epsilon, t}(f)} |Q_\beta (x, y)|(1-|y|)^{\beta - t}dy \right)^p
(1-|x|)^\alpha dx = +\infty.$$
Since $(1-|x|)^t |f_1(x)| \geq (1-|x|)^t|f(x)| - (1-|x|)^t|f(x) - f_1(x)|$ for every $x \in \mathbb B$ we conclude
that $(1-|x|)^t|f_1(x)| \geq (1-|x|)^t |f(x)| - \epsilon_1$ and therefore
$$(\epsilon - \epsilon_1) \chi_{U_{\epsilon, t}(f)}(x)(1-|x|)^{-t} \leq |f_1(x)|, \qquad x \in \mathbb B.$$
Hence
\begin{align*}
+\infty & = \int_{\mathbb B} \left( \int_{U_{\epsilon, t}(f)}|Q_\beta (x, y)|(1-|y|)^{\beta - t} dy \right)^p
(1-|x|)^\alpha dx\\
& = \int_{\mathbb B} \left( \int_{\mathbb B} \frac{\chi_{U_{\epsilon, t}(f)}(y)}{(1-|y|)^t}
|Q_\beta(x, y)|(1-|y|)^\beta dy \right)^p (1-|x|)^\alpha dx \\
& \leq C_{\epsilon, \epsilon_1} \int_{\mathbb B} \left( \int_{\mathbb B} |f_1(y)| |Q_\beta(x, y)|(1-|y|)^\beta dy \right)^p (1-|x|)^\alpha dx = M,
\end{align*}
and we are going to prove that $M$ is finite, arriving at a contradiction. Let $q$ be the exponent
conjugate to $p$. We have, using Lemma \ref{qbeta},
\begin{align*}
I(x) & = \left( \int_{\mathbb B} |f_1(y)| (1-|y|)^\beta |Q_\beta(x, y)| dy \right)^p \\
& = \left( \int_{\mathbb B} |f_1(y)| (1-|y|)^\beta |Q_\beta(x, y)|^{\frac{1}{n+\beta}(\frac{n}{p} + \beta -
\epsilon)} |Q_\beta(x, y)|^{\frac{1}{n+\beta}( \frac{n}{q} + \epsilon)} dy \right)^p \\
& \leq \int_{\mathbb B} |f_1(y)|^p (1-|y|)^{p\beta} |Q_\beta(x, y)|^{\frac{n + p\beta - p\epsilon}{n+\beta}} dy
\left( \int_{\mathbb B} |Q_\beta(x, y)|^{\frac{n+q\epsilon}{n+\beta}} dy \right)^{p/q}\\
& \leq C (1-|x|)^{-p\epsilon}
\int_{\mathbb B} |f_1(y)|^p (1-|y|)^{p\beta} |Q_\beta(x, y)|^{\frac{n + p\beta - p\epsilon}{n+\beta}} dy
\end{align*}
for every $\epsilon > 0$. Choosing $\epsilon > 0$ such that $\alpha - p \epsilon > -1$ we have, by Fubini's theorem
and Lemma \ref{qbeta}:
\begin{align*}
M & \leq C \int_{\mathbb B} |f_1(y)|^p (1-|y|)^{p\beta} \int_{\mathbb B} (1-|x|)^{\alpha - p\epsilon}
|Q_\beta(x, y)|^{\frac{n+p\beta - p\epsilon}{n+\beta}} dx dy \\
& \leq C \int_{\mathbb B} |f_1(y)|^p (1-|y|)^\alpha dy < \infty.
\end{align*}
In order to prove the remaining estimate $t_1(f) \leq C t_2(f)$ we fix $\epsilon > 0$ such that the integral appearing in the definition of $t_2(f)$ is finite and use Theorem on integral representation in unit ball formualted above, with $\beta > \max (t-1, 0)$:
\begin{align*}
f(x) & = \int_{\mathbb B \setminus U_{\epsilon, t}(f)} Q_\beta(x,y) f(y) (1-|y|^2)^\beta dy +
\int_{U_{\epsilon, t}(f)} Q_\beta(x, y) f(y) (1-|y|^2)^\beta dy \\
& = f_1(x) + f_2(x).
\end{align*}
Since, by Lemma \ref{qbeta}, $|f_1(x)|  \leq 2^\beta \int_{\mathbb B} |Q_\beta(x, y)| (1-|w|)^{\beta - t} dy  \leq C(1-|x|)^{-t}$ we have $\| f_1 \|_{A^\infty_t} \leq C \epsilon$. Thus it remains to show that $f_2 \in A^p_\alpha$
and this follows from
$$\| f_2 \|^p_{A^p_\alpha} \leq \| f \|_{A^\infty_t}^p \int_{\mathbb B} \left( \int_{U_{\epsilon, t}(f)}
|Q_\beta (x, y)| (1-|y|^2)^{\beta - t}  dy \right)^p (1-|x|)^\alpha dx < \infty. \quad \Box$$

The above theorem has a counterpart in the $\mathbb R^{n+1}_+$ setting. As a preparation for this result we need the
following analogue of Lemma \ref{qbeta}.

\begin{lem}\label{qm}
For $\delta > -1$, $\gamma > n + 1 + \delta$ and $m \in \mathbb N_0$ we have
$$\int_{\mathbb R^{n+1}_+} |Q_m(z, w)|^{\frac{\gamma}{n+m+1}} s^\delta dy ds \leq C t^{\delta - \gamma + n + 1},
\qquad t > 0.$$
\end{lem}

{\it Proof.} Using Fubini's theorem and estimate (\ref{estq}) we obtain
\begin{align*}
I(t) & = \int_{\mathbb R^{n+1}_+} |Q_m(z, w)|^{\frac{\gamma}{n+m+1}} s^\delta dy ds
\leq C \int_0^\infty s^\delta \left( \int_{\mathbb R^n} \frac{dy}{[|y|^2 + (s+t)^2]^\gamma} \right) ds \\
& = C\int_0^\infty s^\delta (s+t)^{n-\gamma} ds  = Ct^{\delta - \gamma + n + 1}. \qquad \Box
\end{align*}

For $\epsilon > 0$, $\lambda > 0$ and
$f \in h(\mathbb R^{n+1}_+)$ we set:
$$V_{\epsilon, \lambda}(f) = \{ (x, t) \in \mathbb R^{n+1}_+ : |f(x, t)|t^\lambda \geq \epsilon \}.$$

\begin{thm}
Let $p > 1$, $\alpha > -1$, $\lambda = \frac{\alpha + n + 1}{p}$, $m \in \mathbb N_0$ and
Set,
for $f \in \tilde A^\infty_{\frac{\alpha + n + 1}{p}}(\mathbb R^{n+1}_+)$:
$$s_1(f) = {\rm dist}_{\tilde A^\infty_{\frac{\alpha + n +1}{p}}} (f, \tilde A^p_\alpha),$$
$$s_2(f) = \inf \left\{ \epsilon > 0 : \int_{\mathbb R^{n+1}_+} \left( \int_{V_{\epsilon, \lambda}}
Q_m(z, w) s^{m - \lambda}dy ds \right)^p t^\alpha dx dt < \infty \right\}.$$
Then there is an $m_{0}$ depending on $\alpha,p,n$ so that for all $ m > m_{0} $ we have $$s_1(f) \asymp s_2(f)$$.
\end{thm}

The proof of this theorem closely parallels the proof of the previous one, in fact, the role of Lemma \ref{qbeta}
is taken by Lemma \ref{qm} and the role of Theorem on integral representation in unit ball formulated above 
is taken by Theorem on integral reproesentation in upperhalspace which was also given above. We leave details
to the reader.

The following theorem covers the remaining case of $p\leq 1$.
A sharp theorem for this values of $p$ as above in parallel case of same values of $p$, but in unit ball is proved via direct application of Whitney decomposition of upperhalfspace (see \cite{AS2})  we mentioned above 
and estimates for it which were also mentioned above by us (see also\cite{AV}).We remark here that this theorem below is not new and we add this assertion for completness of our exposition
in these issues connected with extremal problems .
The proof of theorem will be omitted. We refer the reader to our previous papers \cite{AS1},\cite{AS2} where the proof of this result can be seen.

\begin{thm}
Let $p \leq 1$, $\alpha > -1$, $\lambda = \frac{\alpha + n + 1}{p}$, $m \in \mathbb N_0$ and
Set,
for $f \in \tilde A^\infty_{\frac{\alpha + n + 1}{p}}(\mathbb R^{n+1}_+)$:
$$s_1(f) = {\rm dist}_{\tilde A^\infty_{\frac{\alpha + n +1}{p}}} (f, \tilde A^p_\alpha),$$
$$s_2(f) = \inf \left\{ \epsilon > 0 : \int_{\mathbb R^{n+1}_+} \left( \int_{V_{\epsilon, \lambda}}
Q_m(z, w) s^{m - \lambda}dy ds \right)^p t^\alpha dx dt < \infty \right\}.$$
Then there is an $m_{0}$ depending from $\alpha,p,n$ so that for all $m>m_{0}$ we have $$s_1(f) \asymp s_2(f)$$.
\end{thm}

The proof closely parallels the proof we have in the unit bàll case  for same values of $p$ (see \cite{AS2})
At the end of this paper we add important remarks on distance theorems in so-called mixed norm and general weighted spaces of harmonic functions.

By $S$ we define a class of all positive measurable functions $v$ on (0,1) for which there are positive numbers $m_{v}$,$M_{v}$ and $q_{v}$ so that
$m_{v},q_{v} \in (0,1)$ and so that $m_{v}\leq \frac {v(\lambda r)}{v(r)}\leq M_{v}$ for all $r$ and $\lambda$,$r\in(0,1)$,$\lambda\in(q_{v},1)$
These spaces of special functions were mentioned in \cite{DS} and analytic Bergman spaces with these general weights were considered and studied in unit disk,polydisk and later in unit ball. \cite{MZ}
This $S$ class includes various unusal weight- functions ,for example, like $t^{\alpha}(\ln(\frac{c}{t})^{\beta}$ for any  positive numbers $\alpha$ and $\beta$.

We define as in \cite{MZ} the following general spaces of harmonic functions of Bergman -type in unit ball and upperhalfspace.
Let $h^{p}_{v}(\mathbb B)$ for all positive values of $p$ be the space of all harmonic $F$  functions in the unit ball with the following finite quazinorm.
$$\int_{\mathbb B}|F(x)|^{p}v(1-|x|)dx$$.These are Banach spaces for all $p$ so that $p>1$ and  quazinormed spaces for all $p\leq 1$ (see \cite{MZ})
We modify this quazinorm in a standard manner to define also $h^{\infty}_{v}(\mathbb B)$ as a space of $F$ harmonic functions with finite quazinorm $$ \sup_{x\in \mathbb B}|F(x)|v(1-|x|)$$

Some results of this note can be extended directly to these general spaces since the Bergman representation formula for these general classes is also valid  (see\cite{MZ}) 
and it is a core of our approach in all proofs \cite{AS1},\cite{AS2},\cite{AS3}.
We formulate a result from \cite{MZ} it says that
for all $p>1$  and $p=\infty$ and all $\alpha> s(v,p)= \frac{{\alpha_{v}}+1}{p} $ ,where
$\alpha_{v}=\frac{log m_{v}}{\ln q_{v}}$ the intergal representation for all functions from $h^{p}_{v}$ via $Q_\alpha(x,y)$ kernel  is valid.
The following theorem in particular can be found in \cite{MZ}.We formulate it  since it is not available in literature openly.

\begin{thm}[\cite{MZ}]\label{intrep}
Let $p \geq 1$ or $p=\infty$ and $\alpha > s(v,p)>0 $. Then for every $f \in h^p_{v}$ and $x \in \mathbb B$ we have
$$f(x) = \int_0^1 \int_{\mathbb S^{n-1}} Q_\alpha(x, y) f(\rho y')(1-\rho^2)^\alpha \rho^{n-1}d\rho dy',
\qquad y = \rho y'.$$
\end{thm}

The complete analogue of this representation is valid for all $p>1$ for spaces of Bergman type in upper halfspace $R^{n+1}$, which we defined above ,but with general $v$ weights from S class \cite{MZ}
These are spaces of $F$ harmonic functions  in upperhalfplane $H^{p}_{v}(R^{n+1}_{+})$ with finite quazinorm $$\int_{R^{n+1}_{+}}|F(x",x_{n+1})|^{p}v(x_{n+1})dxdx_{n+1}$$ for all positive $p$ \cite{MZ}
with obvious modification for $p=\infty$ case.And again for $p>1$ these are Banach spaces, for all other positive $p$ they are quazinormed spaces. 
Here again $v$ is a positive slowly varying function on $(0,\infty)$ from S class.(see \cite {MZ})

In this case another integral repesentation  for upperhalfspace $R^{n+1}_{+}$ with another kernel ,which we also mentioned above is true.
Note for other values of positive $p$,namely for $p\leq 1$ there is also another(more complicated) kind of integral representations \cite{MZ} both  unit ball and upperhalfspace. 
The following result for general weighted Bergman spaces also can be found in \cite{MZ}

\begin{thm}\label{brthm}
Let $1<p \leq\infty$  and  $m>m_{0}$  for certain fixed  $m_{0}$ depending from $p,w$,$m_{0}=s(v,p)-1$
then the following rpresentation is valid
\begin{equation}\label{brep}
f(z) = \int_{\mathbb R^{n+1}_+} f(w) Q_m(z, w) s^m dy ds, \qquad f \in  H^p_{v},\quad z \in \mathbb R^{n+1}_+.
\end{equation}
\end{thm}

Finnaly we remark results of these note partially can be extended to so- called mixed norm spaces of harmonic functions $F^{p,q}_\alpha $
and $B^{p,q}_\alpha $ in unit ball and  their direct analogues in upperhalf space $R^{n+1}_{+}$ at the same time.
For definition of these four scales of  harmonic mixed norm spaces we refer the reader to \cite{AS2},\cite{AS1},\cite{AS3}
Note CArleson -type embedding theorems for these  spaces were  studied  in our previous paper \cite {AS2} . The core of the proof of a distance theorem for these classes  is the same and at the end of proof an appropriate (known)
embeddings connecting classical Bergman $A^{p}_\alpha$  spaces  with $F^{p,q}_{\alpha}$ and $B^{p,q}_{\alpha}$ classes in unit ball or $R^{n+1}_{+}$ should be used \cite{AS2}
We provide two natural examples in $R^{n+1}_{+}$ from a group of not sharp results for these calsses in ball and upperhalfspace
to readers noting we didn"t get yet any sharp result yet in this direction.Note the problem is motivated by an embedding $F^{p,q}_\alpha \subset A^\infty_\lambda $
and also by another embedding $B^{p,q}_\alpha \subset A^\infty_\lambda$ for $q\leq p$ ,which we actually already showed partially above (see also \cite {AS1},\cite{AS2},\cite{AS3}) 
where $ \lambda=\frac{\alpha+n+1}{p}$ ,$\alpha>-1$,$p,q\in(0,\infty)$

\begin{thm}
Let $p,q\in (0,\infty)$,$q \leq p$, $\alpha > -1$, $\lambda = \frac{\alpha + n + 1}{p}$, $m \in \mathbb N_0$ .
Set,for $f \in \tilde A^\infty_{\frac{\alpha + n + 1}{p}}(\mathbb R^{n+1}_+)$:
$$s_1(f) = {\rm dist}_{\tilde A^\infty_{\frac{\alpha + n +1}{p}}} (f, \tilde B^{p,q}_\alpha),$$
$$s_2(f) = \inf \left\{ \epsilon > 0 : \int_{\mathbb R^{n+1}_+} \left( \int_{V_{\epsilon, \lambda}}
Q_m(z, w) s^{m - \lambda}dy ds \right)^p t^\alpha dx dt < \infty \right\}.$$
Then there is an $m_{0}$ depending from $p,q,n,\alpha$ so that for all $m>m_{0}$ ,
$$s_2(f) \leq  C s_1(f)$$.
\end{thm}

\begin{thm}
Let $ q,p \in (0,\infty)$,$q\leq p $, $\alpha > -1$, $\lambda = \frac{\alpha + n + 1}{p}$, $m \in \mathbb N_0$ .
 Set,
for $f \in \tilde A^\infty_{\frac{\alpha + n + 1}{p}}(\mathbb R^{n+1}_+)$:
$$s_1(f) = {\rm dist}_{\tilde A^\infty_{\frac{\alpha + n +1}{p}}} (f, \tilde F^{p,q}_\alpha),$$
$$s_2(f) = \inf \left\{ \epsilon > 0 : \int_{\mathbb R^{n+1}_+} \left( \int_{V_{\epsilon, \lambda}}
Q_m(z, w) s^{m - \lambda}dy ds \right)^p t^\alpha dx dt < \infty \right\}$$
Then there is an $m_{0}$ depending from $p,q,n,\alpha$ so that for all $ m>m_{0}$,
$$s_2(f) \leq C s_1(f)$$
\end{thm}

Proofs repeat arguments we provided above in combination with (known) embeddings between mixed norm and Bergman spaces \cite{AS2} ,namely we use the embedding
$X\subset A^{p}_\alpha$ where by $X$ we denote one of these mixed norm spaces and which we should use at very last step of proof and we omit easy details here.

\end{document}